\documentclass[reqno,10pt]{article}
\usepackage{hyperref}
\usepackage{amsmath}
\usepackage{amssymb}

\newtheorem{thm}{Theorem}[section]
\newtheorem{lemma}[thm]{Lemma}
\newtheorem{corol}[thm]{Corollary}
\newtheorem{propos}[thm]{Proposition}
\newtheorem{rema}{Remark}[section]
\def\bp{\begin{propos}}
\def\ep{\end{propos}}
\def\bt{\begin{thm}}
\def\et{\end{thm}}
\def\bco{\begin{corol}}
\def\eco{\end{corol}}
\def\bl{\begin{lemma}}
\def\el{\end{lemma}}
\def\br{\begin{rema}}
\def\er{\end{rema}}
\def\be{\begin{equation}}
\def\ee{\end{equation}}
\def\ba{\begin{array}}
\def\ea{\end{array}}

\def\P{{\mathbb P}}
\def\E{{\mathbb E}}
\def\R{{\mathbb R}}
\def\Z{{\mathbb Z}}

\def\fA{{\cal A}}

\def\fK{{\cal K}}
\def\fL{{\cal L}}
\def\ze{{\zeta}}

\def\Ga{{\Gamma}}
\def\a{{\alpha}}

\def\QED{\hfill$\square$\vskip 3mm}

\def\Dp{\displaystyle}
\def\hb{\hbox}

\begin{document}

\title{\Large A GEOMETRICAL STRUCTURE FOR AN INFINITE\\  ORIENTED CLUSTER AND ITS
UNIQUENESS\\[10mm]
}

\author{{ Xian-Yuan Wu}\thanks{Research supported in part
by the Natural Science Foundation of China in grant No. 10301023 and
the foundation of ministry of education for experts who have ever studied
abroad.}\and {Yu Zhang}\thanks{Research supported in part by NSF
grant DMS-0405150.}}

\date{}
\maketitle

\begin{abstract}
We consider the supercritical oriented percolation model. Let ${\fK}$ be
all the percolation points. For each $u\in {\fK}$, we write $\gamma_u$ as its right-most path.
Let $G=\cup_u \gamma_u$. In this paper, we show that
 $G$ is a single tree with only one topological end. We also present a relationship between ${\fK}$ and $G$ and
construct a bijection between ${\fK}$ and $\Z$
using the preorder traversal algorithm. Through applications of this fundamental graph property,
we show  the uniqueness of an infinite oriented cluster by ignoring finite vertices.

%To construct the graph
%$G$, we consider all the {\it right most} infinite open pathes
%from points of ${\fK}$, and let $G$ be the random oriented graph
%consisting of all such right most infinite pathes.
\end{abstract}

\section{Introduction and statement of the results.}
\renewcommand{\theequation}{1.\arabic{equation}}
\setcounter{equation}{0}

We consider the graph with vertices ${\fL}=\{(m,n)\in\Z^2:m+n\hbox{
is even}\}$ and oriented edges from $(m,n)$ to $(m+1,n+1)$ and to
$(m-1,n+1)$. The oriented edge from $u$ to $v$ is denoted by
$e(u,v)$.
% and the origin of the lattice $(0,0)$ is denoted by ${\bf 0}$.
As usual, each edge is independently open or closed with a
probability of $p$ or $1-p$. We denote by $\P_p$  the corresponding
product measure and by $\E_p$ the expectation with respect to
$\P_p$. For two vertices $u,v\in {\fL}$, we say $v$ can be reached
from $u$, denoted by $u\rightarrow v$, if there is a sequence of
vertices and open edges $v_0=u,e_1,v_1,\ldots,v_{m-1},e_m,v_m=v$
such that $e_i=e_i(v_{i-1},v_i)$ is open for $1\leq i\leq m$. If there is
no such sequence, we say $v$ cannot be reached from $u$ and denote it
by $u\nrightarrow v$. We define the oriented percolation cluster
at $(x,y)\in {\fL}$ by
$$C_{(x,y)}=\{(z,w)\in {\fL}:(x,y)\rightarrow (z,w)\}.$$
Let
$$\Omega_{(x,y)}=\{|C_{(x,y)}|=\infty\}.
$$ The percolation probability and the critical point are defined by
$$\theta(p)=\P_p(\Omega_{(0,0)})\mbox{ and } \vec {p}_c=\sup\{p:\theta(p)=0\}.$$
It is well known
that $$0<\vec{\hskip0.5mm p}_c<1.$$
By  definition, we know that
$$ \theta(p) =0 \mbox{ if } p< \vec{p}_c\mbox{ and } \theta(p)> 0 \mbox{ if } p> \vec{p}_c.$$
Furthermore, Bezuidenhout and Grimmett  \cite{7} showed that
\be \label{1.1}\theta(p) =0 \Leftrightarrow p\leq \vec{p}_c .\ee

%The above is a version of the general oriented percolation on the
%square lattice $\Z^2$, the notation developed in this paper is
%generally consistent with that of \cite{1,2,3}.

We say $(x,y)\in {\fL}$ is a {\em percolation point} if $|C_{(x,y)}|=\infty$.
By  (\ref{1.1}), if  $p\leq \vec{p}_c$, there is no percolation point,
but if $p > \vec{p}_c$, there are infinitely many percolation points.
When $p>\vec{p}_c$, we collect all
percolation points  and denote them by
\be\label{1.2}{\fK}=\{(x,y)\in{\fL}:
|C_{(x,y)}|=\infty\}.\ee

To understand the  oriented clusters, we need to establish their boundaries.
For $(x,y)\in{\fK}$, let
$\gamma_{(x,y)}$ be the {\it right-most} infinite open path starting
at $(x,y)$. More precisely, let $\gamma_{(x,y)}$ be the infinite
sequence of vertices and open oriented edges
$v_0=(x,y),e_1,v_1,\ldots,e_n,v_n,\ldots$, with $v_n=(x_n,y_n)$ and
$e_{n}(v_{n-1},v_n)$ satisfying
$$\{k\geq 1: (x,y)\rightarrow (x_n+k,y_n)\in {\fK}\}=\emptyset,\mbox{ for all } n\geq
1.$$
Similarly, we may define the {\em left-most} infinite open
path as $\ell_{(x,y)}$, starting at $(x,y)$ by changing $k\geq 1$
in the above equation to $k\leq -1$. With these definitions, note that any
infinite oriented path of $C_{(x,y)}$ will stay in the cone
between $\gamma_{(x,y)}$ and $\ell_{(x,y)}$.

For any infinite oriented open path $\Ga$, and vertex $v\in\Ga$,
let \be\label{brl}\ba{lll}&b_r(v,\Ga):\hskip -3mm
&=\{u\in\fL\setminus{\Ga}:u \mbox{ lies to
the right of }\Ga \mbox{ and }v\rightarrow u \mbox{ uses no edges of }\Ga \},\\[2mm]
& b_l(v,\Ga):\hskip -3mm &=\{u\in\fL\setminus{\Ga}:u \mbox{ lies
to the left of }\Ga \mbox{ and }v\rightarrow u \mbox{ uses no
edges of }\Ga \}.\ea\ee
These will be the right and left {\it buds} of $\Ga$
planted in $v$ (see Fig.~\ref{f2}). For two vertices $u$ and $v$
of $\Ga$ such that $u\rightarrow v$ in $\Ga$, we write $\Ga(u,v)$
for the finite piece of $\Ga$ from $u$ to $v$, and let
\be\label{cr}C_r(\Ga(u,v)):=\Dp\bigcup_{v'\in \Ga(u,v)\setminus
\{v\}}b_r(v',\Ga),\ \ \ C_l(\Ga(u,v)):=\Dp\bigcup_{v'\in
\Ga(u,v)\setminus \{v\}}b_l(v',\Ga).\ee Clearly, if $\Ga$ is a
right-most (resp., left-most) path, then all right (resp., left) buds of
$\Ga$ are finite, which (resp., left-most) implies that $C_r(\Ga(u,v))$ (resp.,
$C_l(\Ga(u,v))$) is finite.

Let $G$ be the random oriented graph consisting of
$\gamma_{(x,y)}$ for $(x,y)\in {\fK}$. In other words, the vertex
set of $G$ is in ${\fK}$ and the edges of  $G$ are  open edges in
$\gamma_{(x,y)}$ for $(x,y)\in {\fK}$. Clearly, there are no loops
in $G$, and $G$ is a {\it forest}.

Each vertex  $u\in {\fL}$ is adjacent to  two edges above $u$, denoted by
{\em upper edges}, and two edges below $u$, denoted by {\em lower edges}.
Note that $G$ consists of  oriented paths without  loops, so each vertex $u$
of $G$ is adjacent to  only one  upper edge in $G$. We call the other endpoint of the upper edge
 the {\em mother vertex} of $u$.
On the other hand, $u$ is also adjacent to one or two  lower  edges in $G$.
%For any $u\in \fK$, ({\bf Is this $G$? not ${\fK}$ })
We   call the other vertex or vertices of the lower edge or
edges of $u$ the {\em  daughter vertex} or the {\em daughter vertices}. By the definition of $G$, every vertex
of $G$ has a mother vertex and at most two daughter vertices. If $u$ has two daughter vertices,
 they are {\it sisters}, and the vertex at the left lower edge
is the {\em older sister} and the vertex at the right lower edge  is the {\em
younger sister}.

Let $M(u)$ denote the mother vertex of $u$, and
iteratively for $n\geq 1$, let $M^n(u)=M(M^{n-1}(u))$ denote
the $n$th ancestor of $u$, where $M^0(u)=u$. We also denote this by
$$
D^n(u,G):=\{v\in{\fK}:M^n(v)=u\} \mbox{ for } n\geq 0 \mbox{ and }
D(u,G):=\Dp\bigcup_{n\geq 0}D^n(u,G)$$
 the $n$th
generation and the set of all descendants of $u$. We call $D(u,G)$
the {\it branch} of $u$.  If
$D(u,G)$ is finite for all $u\in\fK$, we say that $G$ has finite branches. We say that two vertices
$u,v\in \fK$ are connected if they have a common ancestor,
that is, there exist nonnegative integers $n$ and $m$ such that
$M^n(u)=M^m(v)$. This defines an equivalence relation in $\fK$, and
the equivalence classes are called {\it connected components}.
Obviously, a connected component of $G$ is a {\it single tree}.

In graph $G$, a 1-way infinite path is called a {\it ray}. If we
remove  finitely many  vertices from a ray, the rest of the connected part is
still a ray that is called the {\it tail} of the ray. Two rays $R$
and $R'$ are {\it equivalent} if they have the same tail. This is
an equivalence relation on the set of rays in $G$, and the
equivalence classes are called {\it topological ends} (or,
equivalently, graph-theoretical ends) of $G$. Now, with these
definitions, we state a fundamental property for graph $G$ as
follows.

\bt\label{th1.1} For any $p\in(\vec{\hskip0.5mm p}_c,1)$, and  the oriented graph $G$ defined above,  we have
\begin{description}
    \item[\hskip 8mm {\it (i)}] $G$ has a unique connected component,
    \item[\hskip 8mm {\it (ii)}] $G$ has finite branches, and
    \item[\hskip 8mm {\it (iii)}] Each vertex of $G$ has an
    ancestor with a younger sister almost surely.
\end{description}\et

Obviously, items (i) and (ii) of Theorem~\ref{th1.1} tell us
that any two rays $R$ and $R'$ in $G$ have the same tail. This gives the
following corollary.

\bco\label{co1} $G$ has one topological end.\eco

%By Theorem~\ref{th1.1}, we can obtain ({\bf OK?---I think $\fK$ is
%the most essential thing for the model, we get $\fK$ when we know
%the states of all edge of $\fL$, so, I don't think it necessary to
%get $\fK$ from $G$! On the contrary, it's natural to get
%(construct or obtain) $G$ from $\fK$! ----Maybe I have not yet got
%your real meaning!}) $\fK$ in a deterministic way by using $G$.
%{\it unique infinite succession line} as follows.
For any vertex $u$ of $G$, let $\sigma(u)=1$ if $u$ is older than
her sister, and let $\sigma(u)=2$ otherwise. We associate to each
vertex $u$ the sequence of relative sister-order of its ancestors:
let $\sigma_i(u):=\sigma(M^i(u)), i\geq 0$. For any vertices $u$
and $v$ of $G$, because $G$ is a single tree, $u$ and $v$ have
common ancestors. Let $z=M^i(u)=M^j(v)$ be their closest ancestor.
We say $u$ {\it precedes} $v$ if
$\sigma_{i-1}(u)<\sigma_{j-1}(v)$, where
 $\sigma_{-1}(u)=0$.

On the other hand, we define the {\it successor} of $u\in {\fK}$ as $u'\in \fK$ if
$u$ precedes $u'$ and  no vertex precedes  $u$ and
 $u'$. Conversely, $u$ is the {\it predecessor} of $u'$ if
and only if $u'$ is the successor of $u$.

The successor of a vertex can be found by using the following
algorithm. If the vertex has a daughter, we choose the older one.
Otherwise, we move up the tree until we hit the first vertex that
has a younger sister; this younger sister will be the
successor. Note that the existence of such a vertex is guaranteed by  Theorem~\ref{th1.1} (iii).

The predecessor vertex can be also found by using this algorithm.
If the vertex is  the older  sister, her mother will
be the predecessor vertex. If the vertex has an older sister, we
move from her older sister down the tree and choose the younger
daughter at each step until we come to a vertex with no daughter.
This will be the predecessor. Because $G$ has finite branches, the
predecessor vertex can always be found.

We say that there is a {\it succession line} from $u$ to $v$ if
there exists a finite sequence of vertices $u=v_0,\ldots, v_k=v$
such that $v_{i-1}$ is the successor of $v_i$ for $i=1,2,\ldots,k$. We
say $G$ has a {\it unique infinite succession line} if, for every
couple of vertices $u$ and $v$, there is a unique succession line
either from $u$ to $v$ or from $v$ to $u$.

With these definitions, by Theorem~\ref{th1.1}, we have the following
corollary.

\bco\label{co2} $G$ has a unique infinite succession line. The maps
$\Pi(u,G)=u'$ (the successor of $u$) and $\Pi^{-1}(u',G)=u$ (the
predecessor of $u'$) are well-defined, and one is the inverse of
the other. Furthermore, ${\fK}=\{\Pi^n(u,G):n\in\Z\}$ for all
$u\in \fK$.\eco

The succession line was first studied in \cite{10} for {\it Poisson
trees} defined from the two or three-dimensional Poisson point processes.
By constructing succession lines of Poisson trees, Ferrari, Landim and Thorisson \cite{10} proved
the {\it point-stationary} property of the Palm version of the
Poisson point processes. The concept {\em point-stationary} is
defined in \cite{5} and is shown to be the characterizing property
of the Palm version of any stationary point process in $\R^d$.

In view of point processes, $\fK$ is a discrete version of
the translation invariant point process of $\R^2$. On event
$\Omega_{(0,0)}$, if let ${\fK}_0:=\{(x,y)\in\fL:
|C_{(x,y)}|=|C_{(0,0)}|=\infty\}$, then by Theorem 5.1 of \cite{10}
and Corollary~\ref{co2}, we know that ${\fK}_0$ is
point-stationary and ${\fK}_0$ is the Palm version of $\fK$.

With Theorem~\ref{th1.1} in hand, we may try to ask the
question: when properly rescaled, does $G$ converge to the
Brownian Web? (See the definition of the Brownian Web and the related
theorems by Fontes, Isopi, Newman and Ravishankar in \cite{8} and
\cite{9}). Note that by using the convergence criteria given in
\cite{8} and \cite{9}, Ferrari, Fontes, and Wu  \cite{10}
proved that the two-dimensional Poisson trees converge to the
Brownian Web. In fact, one of the original motivations for this
paper was to investigate the convergence of the percolation system
to the Brownian Web.  At this point, we are unable to show this
argument.

By applying this fundamental graph property of $G$, we will
try to characterize  infinite oriented clusters. In fact,  one
of the most important questions  in percolation models is to
investigate the uniqueness of infinite clusters. For  two
different percolation points $(x_1, y_1)$ and $(x_2, y_2)$, Wu
\cite{4} worked on the first step of uniqueness to show
that, for some $(x_3,y_3)\in\fK$,
$$C_{(x_1, y_1)}\cap C_{(x_2, y_2)}\supset C_{(x_3, y_3)}.$$
Recall that for  a general percolation model, Aizeman, Kesten and
Newman (1987) \cite{0} showed the uniqueness of infinite
clusters. We may ask  the uniqueness of infinite oriented clusters.
Clearly, $C_{(x_1, y_1)}\neq C_{(x_2, y_2)}$ for two percolation
points $(x_1, y_1)$ and $(x_2, y_2)$, since we are investigating  oriented
paths. However, even though two infinite oriented clusters are always
different, they might be  different only
 in finitely many vertices. In other words, the main infinite parts of two oriented clusters
are the same. With this observation, we may modify the
definition of  uniqueness to investigate  infinite parts of $C_{(x_1, y_1)}$ and $ C_{(x_2, y_2)}$.
More precisely, for two percolation points $(x_1, y_1)$ and $(x_2, y_2)$, we say
$$C_{(x_1, y_1)}\asymp C_{(x_2, y_2)} \mbox{ if  } C_{(x_1, y_1)}=C_{(x_2, y_2)}\mbox{ except
finitely many vertices.}$$ With this new ``$\asymp$," we may ask what
the uniqueness of infinite oriented clusters is. As expected, we show
the following result.

\bt\label{th1.4} Under the definition of ``$\asymp$," there is only one infinite oriented cluster.
\et

\section{Kuczek's construction, coalescing random walks and CLT.}
\renewcommand{\theequation}{2.\arabic{equation}}
\setcounter{equation}{0}

We use the notations in \cite{2} in Section 2. For $A\subset
(-\infty,\infty)$, we denote a random subset by
$$\xi^A_n=\{x:\exists\ x'\in A \hbox{ such that }(x',0)\rightarrow(x,n)\},\ \
n>0.$$
The right edge of $\xi^{(-\infty,0]}_n$ is defined by
$$r_n=\sup\xi^{(-\infty,0]}_n. \ \ (\hb{Where }\sup\emptyset=-\infty.) $$

We know (see page 1004 in \cite{1}) by using a subadditive argument that there exists a nonrandom constant $\alpha(p)$
such that
$$\lim_{n\rightarrow\infty}\frac{r_n}n=\inf_n\left\{\frac{\E_p(r_n)}{n}\right\}=\alpha(p)\hbox{ a.s. and in }L_1.$$
It has been proved in \cite{7,1} that
$$\alpha(p)=-\infty,\hbox{ if }p<\vec{\hskip0.5mm p}_c,\hbox{ and
}\alpha(\vec{\hskip0.5mm p}_c)=0,
\hbox{ and }1\geq\alpha(p)>0\hbox{ if }p>\vec{\hskip0.5mm p}_c.$$
In particular,  $\alpha(p)$ is infinitely differentiable for
all $p\in (\vec{\hskip0.5mm p}_c, 1)$ (see \cite{3}).

Let us now denote
$$ \xi'_0=\xi^{\{(0,0)\}}_0,$$ and for all $n\geq 0$,
$$\xi'_{n+1}=\left\{\ba{lll}&\hskip-5mm\{x:(y,n)\rightarrow(x,n+1)\hb{ for some }
y\in \xi'_n\},&\hbox{ if this set is non-empty;}\\[3mm]&\hskip-5mm\{n+1\}, &
\hbox{ if otherwise.}\ea\right.$$
Let $r'_n=\sup\xi'_n$. On event $\Omega_{(0,0)}$,
%$\{\xi_n^{\{(0,0)\}}\not=\emptyset\}$,
we know that
$$r'_n=r_n=\gamma_{(0,0)}(n),$$
where $\gamma_{(0,0)}(n)= \gamma_{(0,0)}\cap \{y=n\}$.

\vskip 2mm

Let $T_0=0$ and $T_m=\inf\{n\geq T_{m-1}+1:(r'_n,n)\in\fK\}$ for
$m\geq 1$. Define $\tau_0=0$ and $$\tau_1=T_1,
\tau_2=T_2-T_1,\ldots,\tau_m=T_m-T_{m-1},\ldots,$$ where
$\tau_i=0$ if $T_i$ and $T_{i-1}$ are infinity. Also define
$X_0=0$ and $$X_1=r'_{T_1}, X_2=r'_{T_2}-r'_{T_1},\ldots,
X_m=r'_{T_m}-r'_{T_{m-1}},\ldots,$$ where $X_i=0$ if $T_i$ and
$T_{i-1}$ are infinity. The collection of points $(r'_{T_m},T_m),m\geq
0$ are called {\it break points} for  point $(0,0)$.

In the case of $p\in(\vec{\hskip0.5mm p}_c,1)$, with these
definitions, Kuczek  \cite{2} proved  the following
proposition. \bp\label{p2.1}
%{\rm
For $p\in(\vec{\hskip0.5mm p}_c,1)$, on $\Omega_{(0,0)}$,
$\{(X_m,\tau_m):m\geq 1\}$ are independently identically
distributed with all moments, and
%\bp\label{2.2} Suppose that $p\in(\vec{\hskip0.5mm p}_c,1)$, then
%conditioned on $\Omega_{0}$,
$\Dp\frac{\gamma_{(0,0)}(n)-\a(p)n}{\sqrt{n\sigma^2}}$ converges to
$N(0,1)$ in distribution as $n\rightarrow\infty$, where
$\sigma^2=\E(X_1\E\tau_1-\tau_1\E X_1)^2>0$ and $\E$ is the
expectation with respect to the conditional measure
$\P_p(\cdot\mid \Omega_{(0,0)})$.\ep

Now we turn to our right-most infinite paths. On
$\Omega_{(0,0)}$, the right-most infinite path $\gamma_{(0,0)}$ is
well defined and all break points $(r_{T_m},T_m), m\geq 1$ are
well embedded in $\gamma_{(0,0)}$. By Proposition~\ref{p2.1}, on
event $\Omega_{(0,0)}$, we define an integer-valued random walk
$\zeta_{(0,0)}=\{\zeta_{(0,0)}(t):t\geq 0\}$ as follows:
$$\zeta_{(0,0)}(0):=0 \mbox{ and }\zeta_{(0,0)}(t):=\sum_{i=0}^{N(t)}X_i,\ \hbox{for } t>0,$$
where $N(t)$ is
the largest integer $m$ such that $T_m\leq t$. Note that a break
point defined above is also a  jump point of $\ze_{(0,0)}$. In the
same way, we define the random walk
$\zeta_{(x,y)}:=\{\zeta_{(x,y)}(t):t\geq y\}$.  On event
$\Dp\Omega_{(x,y)}$, $(x,y)\in\fL$. For any vertices
$u_1=(x_1,y_1), u_2=(x_2,y_2),\ldots,u_k=(x_k,y_k)$ in $\fL$, on
event $\Omega_{u_1}\cap\ldots\cap\Omega_{u_k}$, we say random
walks $\zeta_{u_{i_1}}, \zeta_{u_{i_2}}$ {\it meet} if  two walks {\it jump
synchronously} to the  same position  at some
time $t_0(\in\Z)\geq y_{i_1}\vee y_{i_2}$. By this definition of meeting,
once two walks meet, they will coalesce into one thenceforth. This
defines a finite system of coalescing random walks. For random
walks and our right-most infinite paths, we have the following
proposition to describe the relationship between the jump points in two random walks and
the meeting points in two right-most open paths.
%We point out that this proposition will not
%be used to show our theorems, but  for independent interests.

%({\bf It seems we do not use this Prop 2.3 in the
%proof of Thm. Am I right? ------Yes! I only use Proposition 2.3 to
%display the relationship between G and a kind of coalescing random
%walk system, I think it be independently interested!})

\bp\label{p2.2} For any $p\in(\vec{\hskip0.5mm p}_c,1)$ and any
pair $u_1, u_2\in \fL$, conditioned on
$\Omega_{u_1}\cap\Omega_{u_2}$, the following two statements (i)
and (ii) are equivalent, where
\begin{description}
    \item[\hskip 8mm{\it (i)}] $\ze_{u_1}$ meets $\ze_{u_2}$,
    \item[\hskip 8mm{\it (ii)}] the right-most infinite paths $\gamma_{u_1}$ and
$\gamma_{u_2}$ meet.
\end{description}
\ep

{\bf Proof.}  It is clear that (i) implies (ii), so it suffices to prove that
(ii) implies (i).

Let $\Ga_1$ and $\Ga_2$ be two realizations of $\gamma_{u_1}$ and
$\gamma_{u_2}$, respectively. Without loss of generality, we may assume
that $\Ga_1$ and $\Ga_2$ meet at $u_{1,2}\in\fK$ and $u_1$ {\it
precedes} $u_2$. Recall that the concept of {\it precedes} is defined
in Section 1.

Note that  $u_{1,2}$ is a break point for $u_1$, so it is also  a
jump point for $\ze_{u_1}$. To prove  Proposition~\ref{p2.2}, we
find a
 point $v_{1,2}\in \Ga_1\cap\Ga_2$ preceding $u_{1,2}$
such that $v_{1,2}$ is a common jump point of $\ze_{u_1}$ and
$\ze_{u_2}$.

Now let $\{v_m=(x_m,y_m):m\geq 0\}$ be jump (or break) points of
$\ze_{u_2}$ such that $v_0=u_2$, $v_1$ precedes $v_0$, $\cdots$,
and $v_{m}$ precedes $v_{m-1}$. If $u_{1,2}$ is also one of the jump
points for $\ze_{u_2}$, then $u_{1,2}$ should be the meeting point
of $\ze_{u_1}$ and $\ze_{u_2}$. Proposition~\ref{p2.2} follows by
taking $v_{1,2}=u_{1,2}$. If $u_{1,2}$ is not a jump point for
$\ze_{u_2}$, then there exists some $k\geq 1$ such that
$u_{1,2}\in\Ga_2(v_{k-1},v_k)$. By the definition of break point,
$y_k$, which is the second coordinate of $v_k$, $>y$ for any $(x,y)\in
C_r(\Ga_2(v_{k-1},v_k))$, where $C_r(\Ga_2(v_{k-1},v_k))$ is
defined in (\ref{cr}). Note that
$C_r(\Ga_1(u_{1,2},v_k))=C_r(\Ga_2(u_{1,2},v_k))\subset
C_r(\Ga_2(v_{k-1},v_k))$, so $y_k>y$ for any $(x,y)\in
C_r(\Ga_1(u_{1,2},v_k))$. By the fact that $u_{1,2}$ is a jump
(or break) point for $\ze_{u_1}$ and by the definition of break point,
we know that $v_k$ is a jump (or break) point of $\ze_{u_1}$.
Proposition~\ref{p2.2} follows by taking $v_{1,2}=v_k$. \QED

Write $R_{\a(p)}$ for the line in $\R^2$ with the equation $\Dp
y={x}/{\a(p)}$. Conditioned on $\Omega_{(0,0)}$, let us consider
the behavior of the right-most infinite path $\gamma_{(0,0)}$. By
Proposition~\ref{p2.1}, we have the following proposition.

\bp\label{p2.3} Suppose that $p\in(\vec{\hskip0.5mm p}_c,1)$,
then, conditioned on $\Omega_{(0,0)}$, almost surely, the
right-most infinite path $\gamma_{(0,0)}$ crosses the line $R
_{\a(p)}$ infinite many times.\ep

%Note that we are working on the planar graph. The main advantage is that the right most
%oriented open is well defined. In other words, on event $\{xi_n^{(0,0)}\neq \emptyset\}$, we know that
%$r_n'$ is well defined.
%The same terminology is also used in general percolation as the lowest open
%crossing. Note that (see Chapter 9 in [Grimmett]) when
%When $r_n=\Gamma_n$ for a fixed oriented path from $(0,0)$ to the horizontal line $y=n$,
%we may convenced  ourselvies with a few diagrams that event $\{r_n=\Gamma_n\}$ only depends
%on the configurations of the edges in  the right of $\Gamma_n$. We state this argument as the following proposition.
%For a fixed oriented path $\Gamma_n$ that starts with $(0,0)$, let $R(\Gamma)$ be all edges in $\fL$ in the
%right side of $\Gamma_n$ and below the horizontal line $y=n$, where we assume that
%$R(\Gamma_n$ also contains the edges in $\Gamma_n$ and $y=n$.

%\bp\label{2.4} The event $\{r_n'=\Gamma_n\}\in \Prod_{e\in (R(\Gamma_n))}$(open, closed) for a fixed oriented path from $(0,0)$ to $y=n$.\ep

\section{Proof of Theorem~\ref{th1.1}.}
\renewcommand{\theequation}{3.\arabic{equation}}
\setcounter{equation}{0}

Before our proofs, we need to introduce a few  notations (see
Fig.~\ref{f1}). For any $u=(x,y)\in \fL$, we define
$$\vee{_{u}}:=\{v\in{\fL}:
\hb{ there is an oriented path from $u$ to $v$}\}.$$
Note that
the oriented path in the definition of $\vee{_{u}}$ does not need
 to be open, so $C_u \subset \vee{_{u}}$. Similarly, we define
$$\wedge{_{u}}:=\{v\in{\fL}:
\hb{ there is an oriented path from $v$ to $u$}\}.$$
For any
$n\geq 0$ and $u=(x,y)$, let
$$\wedge{_{u}}(n):=\{v=(x',y')\in \wedge{_{u}}: y-y'\leq n\}.$$
For a finite set $A\subset\fL$ contained in a horizontal line, we
define
$$\wedge{_{A}}:=\{v\in{\fL}:
\hb{ there is an oriented path from $v$ to some point $u$ of
$A$}\},$$ and
$$\wedge{_{A}}(n):=\{v=(x',y')\in \wedge{_{A}}: y-y'\leq n\},\ \ \hb{for }n\geq
0,$$
where $y$ is the second coordinate of some vertex in $A$.

In our proofs, we need to use the following {\it anti-oriented
open path}. Given $u$ and $v\in \wedge{_{u}}$, we say there is an
anti-oriented open path from $u$ to $v$, if $v\rightarrow u$. For
any $u\in\fL$, let
$$C^{\rm anti}_{u}:=\{v\in{\fL}:v\rightarrow u\};$$ clearly, $C^{\rm
anti}_{u}$ is a random subset of $\wedge{_{u}}$. On event $|C^{\rm
anti}_{u}|=\infty$, we write $\ell^{\rm anti}_u$ for the {\it
left-most} anti-oriented infinite open path from $u$.

When $p>\vec{\hskip 0.5mm p}_c$, we have
\be\label{3.1}\P_p(|C^{\rm anti}_u|=|C_u|=\infty)=\theta(p)^2>0,
\forall\ u\in\fL.\ee Vertex $u$, satisfying $|C^{\rm
anti}_u|=|C_u|=\infty$, is called a {\it bidirectional percolation
point}, denoted by $\tilde{\fK}$, the set of all bidirectional
percolation points.

\vskip 3mm {\it Proof of Theorem~\ref{1.1} (i).} It suffices to
prove that, for any vertices $u_1=(x_1,y_1)$,
$u_2=(x_2,y_2)\in\fL$ with $y_1=y_2$, conditioned on
$\Omega_{u_1}\cap\,\,\Omega_{u_2}$, $\gamma_{u_1}$ and $\gamma_{u_2}$
will meet almost surely. In fact, in the case that $y_1\not=y_2$,
on $\Omega_{u_1}\cap\Omega_{u_2}$, there exists some $u_i'=(x_i',y_i'), i=1,2$ almost surely
in $\fK$ with
$y_1'=y_2'=b$ such that
$$x_1'<\min\{x:(x,b)\in\vee{_{u_1}}\cup\vee{_{u_2}}\}\leq
\max\{x:(x,b)\in\vee{_{u_1}}\cup\vee{_{u_2}}\}<x_2'.$$
If $\gamma_{u_1'}$ and $\gamma_{u_2'}$ meet, then
$\gamma_{u_1}$ and $\gamma_{u_2}$ meet.

By translation invariance, we choose $u_2=(0,0)$, and
$u_1=(-n_0,0)$ for some $n_0\geq1$. Let us focus  on
$\gamma_{(0,0)}$. For any realization of $\Ga$ of $\gamma_{(0,0)}$,
by Proposition~\ref{p2.3}, we may assume that $\Gamma$ crosses the
line $R_{\a(p)}$ infinitely many times. For some vertex
$v\in\Gamma$, let $e(u,v)$ be the lower  (oriented) edge of $v$ in
$\Gamma$. We call $v$ a {\it crossing point} if $e(u,v)\cap
R_{\a(p)}\not=\emptyset$.

Given such a realization of $\Gamma$, we  define a series of
independent events $E(k,\Gamma), k\geq 1$ as follows.

We fix $\epsilon_0>0$ such that
$$\Dp\int_{-\infty}^{-\epsilon_0}\frac1{\sqrt{2\pi}}\exp{(-\frac12x^2)}dx>\frac13.$$
By Proposition~\ref{p2.2}, we choose $N_0$ large enough such that
\be\label{3.2}\P_p(\gamma_{(0,0)}(n^2)-\a(p)n^2<-n\sigma\epsilon_0\mid
\Omega_{(0,0)})\geq\frac 13,\ \ \hb{for all $n\geq N_0$}.\ee Let
$v_0(\Ga)=(x_0,y_0)=(0,0)$. We go along $\Gamma$ from $(0,0)$ to
meet  $v_1(\Gamma)=(x_1,y_1)$ (see Fig.~\ref{f1}), one of crossing
points, with $y_1>\max\{{n_0}/({\epsilon_0\sigma}),N_0\}^2$.
Iteratively, we go along $\Gamma$ from $v_{k-1}(\Ga)$ to meet
$v_k(\Gamma)=(x_k,y_k)$, one of crossing points, with
%({\bf This may cause confusion: $v_k(\Ga)$ is not
%the $k$-th crossing point, but the $k$-th crossing satisfying our
%condition.What do you think about it?}), with
\be\label{3.3}y_k-y_{k-1}>
\max\{{2y_{k-1}}/({\epsilon_0\sigma}),N_0\}^2.\ee
%In the
%definition of $v_k(\Gamma)$, $k\geq1$, one may require $y_k$ as
%small as possible.

\begin{figure}
\unitlength=0.4mm
\begin{picture}(280,100)(-70,0)
\put(-10,0){\line(1,0){270}}\put(80,0){\circle*{2}}
\put(82,-8){(0,0)}\put(70,0){\circle*{2}}\put(65,0){\line(-1,1){22}}
\put(68,-8){$u_2$}\put(80,0){\line(1,3){52}}
\put(80,0){\line(0,1){8}}\put(80,8){\line(1,2){21}}\put(101,50){\line(-1,2){10}}\put(91,70){\line(1,1){29}}
\put(120,99){\line(0,1){28}}\put(120,127){\vector(1,1){20}}\put(130,130){$\Ga$}
\put(65,0){\line(1,1){22}}\put(87,22){\line(1,-1){22}}\put(87,22){\circle*{2}}\put(93,15.5){\small
$v_1(\Ga)$}
\put(87,22){\line(-1,0){80}}\put(87,22){\line(1,0){160}}
\put(0,0){\line(1,1){120}}\put(120,120){\line(1,-1){120}}\put(120,120){\circle*{2}}\put(123,118){$v_2(\Gamma)$}
\put(108,145){$R_{\a(p)}$}
\put(105,55){$\wedge_{{v_2(\Ga)}}(y_2-y_{1})$}

%\multiput(66,1)(2,0){21}{\line(1,0){1}}\multiput(68,3)(2,0){19}{\line(1,0){1}}
%\multiput(70,5)(2,0){17}{\line(1,0){1}}\multiput(72,7)(2,0){15}{\line(1,0){1}}
%\multiput(74,9)(2,0){13}{\line(1,0){1}}\multiput(76,11)(2,0){11}{\line(1,0){1}}
%\multiput(78,13)(2,0){9}{\line(1,0){1}}\multiput(80,15)(2,0){7}{\line(1,0){1}}
%\multiput(82,17)(2,0){5}{\line(1,0){1}}\multiput(84,19)(2,0){3}{\line(1,0){1}}
%\multiput(86,21)(2,0){1}{\line(1,0){1}}

\multiput(66,0)(0,2){1}{\line(0,1){1}}\multiput(68,0)(0,2){2}{\line(0,1){1}}\multiput(70,0)(0,2){3}{\line(0,1){1}}
\multiput(72,0)(0,2){4}{\line(0,1){1}}\multiput(74,0)(0,2){5}{\line(0,1){1}}\multiput(76,0)(0,2){6}{\line(0,1){1}}
\multiput(78,0)(0,2){7}{\line(0,1){1}}\multiput(80,0)(0,2){8}{\line(0,1){1}}\multiput(82,0)(0,2){9}{\line(0,1){1}}
\multiput(84,0)(0,2){10}{\line(0,1){1}}\multiput(86,0)(0,2){11}{\line(0,1){1}}
\multiput(88,0)(0,2){11}{\line(0,1){1}}
\multiput(108,0)(0,2){1}{\line(0,1){1}}\multiput(106,0)(0,2){2}{\line(0,1){1}}
\multiput(104,0)(0,2){3}{\line(0,1){1}}\multiput(102,0)(0,2){4}{\line(0,1){1}}\multiput(100,0)(0,2){5}{\line(0,1){1}}
\multiput(98,0)(0,2){6}{\line(0,1){1}}\multiput(96,0)(0,2){7}{\line(0,1){1}}\multiput(94,0)(0,2){8}{\line(0,1){1}}
\multiput(92,0)(0,2){9}{\line(0,1){1}}\multiput(90,0)(0,2){10}{\line(0,1){1}}

\put(20,8){$\wedge_{A_2}(y_{1})$}
\multiput(1,1)(2,0){32}{\line(1,0){0.5}}\multiput(3,3)(2,0){30}{\line(1,0){0.5}}\multiput(5,5)(2,0){28}{\line(1,0){0.5}}
\multiput(7,7)(2,0){7}{\line(1,0){0.5}}\multiput(49,7)(2,0){5}{\line(1,0){0.5}}
\multiput(9,9)(2,0){6}{\line(1,0){0.5}}\multiput(49,9)(2,0){4}{\line(1,0){0.5}}
\multiput(11,11)(2,0){5}{\line(1,0){0.5}}\multiput(49,11)(2,0){3}{\line(1,0){0.5}}
\multiput(13,13)(2,0){4}{\line(1,0){0.5}}\multiput(49,13)(2,0){2}{\line(1,0){0.5}}
\multiput(15,15)(2,0){18}{\line(1,0){0.5}}\multiput(17,17)(2,0){16}{\line(1,0){0.5}}
\multiput(19,19)(2,0){14}{\line(1,0){0.5}}\multiput(21,21)(2,0){12}{\line(1,0){0.5}}
\multiput(23,23)(2,0){97}{\line(1,0){0.5}}\multiput(25,25)(2,0){95}{\line(1,0){0.5}}
\multiput(27,27)(2,0){93}{\line(1,0){0.5}}\multiput(29,29)(2,0){91}{\line(1,0){0.5}}
\multiput(31,31)(2,0){89}{\line(1,0){0.5}}\multiput(33,33)(2,0){87}{\line(1,0){0.5}}
\multiput(35,35)(2,0){85}{\line(1,0){0.5}}\multiput(37,37)(2,0){83}{\line(1,0){0.5}}
\multiput(39,39)(2,0){81}{\line(1,0){0.5}}\multiput(41,41)(2,0){79}{\line(1,0){0.5}}
\multiput(43,43)(2,0){77}{\line(1,0){0.5}}\multiput(45,45)(2,0){75}{\line(1,0){0.5}}
\multiput(47,47)(2,0){73}{\line(1,0){0.5}}\multiput(49,49)(2,0){71}{\line(1,0){0.5}}
\multiput(51,51)(2,0){69}{\line(1,0){0.5}}
\multiput(53,53)(2,0){26}{\line(1,0){0.5}}\multiput(163,53)(2,0){12}{\line(1,0){0.5}}
\multiput(55,55)(2,0){25}{\line(1,0){0.5}}\multiput(163,55)(2,0){11}{\line(1,0){0.5}}
\multiput(57,57)(2,0){24}{\line(1,0){0.5}}\multiput(163,57)(2,0){10}{\line(1,0){0.5}}
\multiput(59,59)(2,0){23}{\line(1,0){0.5}}\multiput(163,59)(2,0){9}{\line(1,0){0.5}}
\multiput(61,61)(2,0){22}{\line(1,0){0.5}}\multiput(163,61)(2,0){8}{\line(1,0){0.5}}
\multiput(63,63)(2,0){57}{\line(1,0){0.5}}\multiput(65,65)(2,0){55}{\line(1,0){0.5}}
\multiput(67,67)(2,0){53}{\line(1,0){0.5}}\multiput(69,69)(2,0){51}{\line(1,0){0.5}}
\multiput(71,71)(2,0){49}{\line(1,0){0.5}}\multiput(73,73)(2,0){47}{\line(1,0){0.5}}
\multiput(75,75)(2,0){45}{\line(1,0){0.5}}\multiput(77,77)(2,0){43}{\line(1,0){0.5}}
\multiput(79,79)(2,0){41}{\line(1,0){0.5}}\multiput(81,81)(2,0){39}{\line(1,0){0.5}}
\multiput(83,83)(2,0){37}{\line(1,0){0.5}}\multiput(85,85)(2,0){35}{\line(1,0){0.5}}
\multiput(87,87)(2,0){33}{\line(1,0){0.5}}\multiput(89,89)(2,0){31}{\line(1,0){0.5}}
\multiput(91,91)(2,0){29}{\line(1,0){0.5}}\multiput(93,93)(2,0){27}{\line(1,0){0.5}}
\multiput(95,95)(2,0){25}{\line(1,0){0.5}}\multiput(97,97)(2,0){23}{\line(1,0){0.5}}
\multiput(99,99)(2,0){21}{\line(1,0){0.5}}\multiput(101,101)(2,0){19}{\line(1,0){0.5}}
\multiput(103,103)(2,0){17}{\line(1,0){0.5}}\multiput(105,105)(2,0){15}{\line(1,0){0.5}}
\multiput(107,107)(2,0){13}{\line(1,0){0.5}}\multiput(109,109)(2,0){11}{\line(1,0){0.5}}
\multiput(111,111)(2,0){9}{\line(1,0){0.5}}\multiput(113,113)(2,0){7}{\line(1,0){0.5}}
\multiput(115,115)(2,0){5}{\line(1,0){0.5}}\multiput(117,117)(2,0){3}{\line(1,0){0.5}}
\multiput(119,119)(2,0){1}{\line(1,0){0.5}}
\put(118,8){\vector(-1,0){12}}\put(120,5){${\rm Area}(1,\Ga)$}
\put(188,70){\vector(-1,0){12}}\put(190,67){${\rm Area}(2,\Ga)$}
\put(19,35){$A_2$}\put(24,33){\vector(1,-1){10}}\put(48,-12){$A_1$}\put(56,-9){\vector(3,2){10}}
\linethickness{1pt}
\put(70,0){\line(-1,0){5}}\put(43,22){\line(-1,0){21}}\put(43,22){\circle*{2}}
\put(22,22){\circle*{2}}\put(65,0){\circle*{2}}
\end{picture}
\vskip 5mm
\begin{center}
\begin{minipage}{13.5cm}
\caption{\small This figure reveals that ${\rm Area}(1,\Ga)$ and
${\rm Area}(2,\Ga)$ are edge-disjoint areas. The situation for large
$k$ is similar.} \label{f1}
\end{minipage}
\end{center}
\end{figure}
Now we define $E(1,\Ga)$ to be the event that there is an
anti-oriented open path from $v_1(\Ga)$ to the half line
$(-\infty,-n_0]\times\{0\}$. Iteratively, for $k\geq 2$, we define
$E(k,\Ga)$ to be the event that there is an anti-oriented open
path from $v_k(\Ga)$ to the half line
$(-\infty,x_k-2y_{k-1}]\times\{y_{k-1}\}$ and then to the line
$(-\infty,+\infty)\times\{0\}$. By this definition,
$E(k,\Ga),k\geq 1$ depends on the edges (see Fig.~\ref{f1}) in
${\rm Area}(k,\Ga)=\wedge_{{v_k(\Ga)}}(y_k-y_{k-1})\cup
\wedge_{{A_k}}(y_{k-1})$, where
$$\ba{ll}A_1\hskip-3mm&=\{u=(x,y)\in\wedge_{v_1(\Ga)}(y_1):x\leq
-n_0,y=0\}(=\wedge_{A_1}(y_0)=\wedge_{A_1}(0)),\\[2mm]
         A_k\hskip-3mm&=\{u=(x,y)\in\wedge_{{v_k(\Ga)}}(y_k-y_{k-1}):x\leq
x_{k-1}-2y_{k-1},y=y_{k-1}\},\ k\geq 2.\ea
         $$
Note that for Area($k,\Ga$), $k\geq 1$ are edge-disjoint areas, so for
a fixed oriented path $\Gamma$, $E(k,\Ga), k\geq1$ are independent
(see Fig.~\ref{f1}). Furthermore, by (\ref{3.2}), (\ref{3.3}), and
by using Proposition~\ref{p2.2} for $\ell^{\rm anti}_{v_k(\Ga)}$,
the left-most anti-oriented infinite open path from $v_k(\Ga)$, we
have \be\label{3.4}\ba{ll}\P_p(E(k,\Ga))&\hskip-3mm\geq
\P_p(|C^{\rm anti}_{v_k(\Ga)}|=\infty \hb{ and } \ell^{\rm
anti}_{v_k(\Ga)}\cap
A_k\not=\emptyset)\\[3mm]
&\hskip-3mm\geq\P_p(\ell^{\rm anti}_{v_k(\Ga)}\cap
A_k\not=\emptyset \mid |C^{\rm anti}_{v_k(\Ga)}|=\infty
)\cdot\P_p(|C^{\rm anti}_{v_k(\Ga)}|=\infty)\\[3mm]
&\hskip-3mm\geq \Dp \frac13\theta(p)>0.\ea\ee

We point out here that $\Ga$ is only used to determine the vertex
set $\{v_k(\Ga):k\geq 0\}$. Furthermore, we need to work on the event
family $\{E(k,\Ga):k\geq1\}$ on event $\gamma_{(0,0)}=\Ga$.
%({\bf Note that $P(\gamma_{(0,0)}=\Ga)=0$, so we may have trouble
%to decompose by using $\cup_{\Gamma} \{\gamma_{(0,0)}=\Ga\}$. What
%do you think? To get rid this trouble, we may take $\Gamma$ very
%long-----Here we only use condition probability, I think we can do
%like that even though $P(\gamma_{(0,0)}=\Ga)=0$, in fact, we use a
%integration other than a summation. If we use $\Gamma$ very long,
%of course, it is right, but it will make the argument more
%complicated!)}
Now, on event $\{\gamma_{(0,0)}=\Ga\}$, let $E^*(1,\Ga)$ be the
event that there is an anti-oriented open path from
$\Ga(0,v_1(\Ga))$ to $A_1$,  and let $E^*(k,\Ga), k\geq 2$ be the
event that there is an anti-oriented open path from
$\Ga(v_{k-1}(\Ga),v_1(\Ga))$ to $A_k$ and then to
$(-\infty,+\infty)\times \{0\}$. On $\{\gamma_{(0,0)}=\Ga\}$. Note
that $\Ga$ is open, so event $E^*(k,\Ga)$ only depends on the
configurations of the edges of ${\rm Area}(k,\Ga)$ lying on the left
side of $\Ga$ (see Fig.~\ref{f1}). Note also that event
$\{\gamma_{(0,0)}=\Gamma\}$ can be decomposed into the intersection
of the following two events:
\begin{enumerate}
    \item ${\cal A}=\{\Gamma$ is open\};
    \item ${\cal B}=\{v\not \rightarrow\infty \mbox{ in $R(\Gamma)$  for each vertex $v\in
    \Gamma$}\}$,
\end{enumerate}
where $R(\Gamma)$ is the edge set to the right of $\Gamma$. It follows from the definition
that ${\cal B}$ only depends on the configurations of the edges in
$R(\Gamma)$. By this decomposition, we have
%Similar to $E(k,\Ga),k\geq1$,
%on $\gamma_{0}=\Ga$, $E^*(k,\Ga),k\geq1$ are
%independent.
\be\label{3.5}\P_p(E^*(k,\Ga)\mid
\gamma_{0}=\Gamma)=\P_p(E(k,\Ga)\mid {\cal A}).\ee For any $k\geq
1$, let $\fA_k(\Ga)$ be the event that all edges in
$\Ga(v_{k-1}(\Ga),v_k(\Ga))$ are open. Note that $\fA_k(\Ga)$,
$k\geq 1$ are increasing events, so by the FKG inequality,
\be\label{3.6}\P_p(E(k,\Ga)\mid \Gamma\mbox{ is
open})=\P_p(E(k,\Ga)\mid {\cal A})=\P_p(E(k,\Ga)\mid
\fA_k(\Ga))\geq \P_p(E(k,\Ga)).\ee
%Therefore, by (\ref{3.5}) and \be\label{3.6}\P_p(E^*(k,\Ga)\mid
%\gamma_{0}=\Gamma)=\P_p(E(k,\Ga)\mid {\cal A})\geq
%.\ee
%{\bf Here $\fA$ is a zero-measured event, do
%you think it convenient to use some version of FKG inequality to
%get the above? Maybe we can insert a step like $$\P_p(E(k,\Ga)\mid
%{\cal A})=\P_p(E(k,\Ga)\mid {\cal A}(k)), $$ and then use the
%usual FKG inequality for increasing events $E(k,\Ga)$ and ${\cal
%A}(k)$. Where ${\cal A}(k)$ be the event that all edges in
%$\Ga(v_{k-1}(\Ga),v_k(\Ga))$ are open, clearly, it has positive
%measure}).
%({\bf I think it's better to use $\fA_k(\Ga)$ instead of
%$B_k(\Ga)$! })
 On the other hand, by the same argument of (\ref{3.5}),
we know that $E^*(k,\Ga)$ for $k\geq 1$ depend on the different
edge layers to the left of $\Ga$ (see Fig.~\ref{f1}). Therefore,
$E^*(k,\Ga)$ for $k\geq 1$  are independent and have
 probabilities bounded away from $\theta(p)/3$ on $\{\gamma_{(0,0)}=\Ga\}$.
 With these observations, by (\ref{3.4}) and the
Borel-Cantelli second lemma, on $\gamma_{(0,0)}=\Ga$,
$E^*(k,\Ga),k\geq 1$ occur infinitely often almost surely. Using
the definition of right-most open path, we know that if
$E^*(k,\Ga)$ occurs for some $k$, then $\gamma_{u_1}$ will meet
$\gamma_{(0,0)}=\Ga$ in $\Ga(0,v_k(\Ga))$. Thus, this shows that
$\gamma_{u_1}$ and $\gamma_{u_2}$ meet, so Theorem~\ref{th1.1} (i)
follows.\QED

{\it Proof of Theorem~\ref{th1.1} (ii).} For any $u\in\fK$, by the
definition of $D(u,G)$, we know that $D(u,G)\subset C^{\rm
anti}_u$. Note that if $u\in{\fK}\setminus\tilde{\fK}$, that is, $u$
is a percolation point but not a bidirectional percolation point,
then $|C^{\rm anti}_u|<\infty$ and $|D(u,G)|<\infty$, so it
suffices to prove that $|D(u,G)|<\infty$ for $u\in\tilde{\fK}$.

By translation invariance, it suffices to prove that
$|D((0,0),G)|<\infty$ almost surely when $(0,0)\in\tilde{\fK}$.

Let $\ell^{\rm anti}_{(0,0)}$ be the left-most anti-oriented
infinite open path from $(0,0)$ and let $L^{\rm anti}$ be a
possible realization of $\ell_0^{\rm anti}$ crossing $R_{\a(p)}$
infinitely many times. Then, using Proposition~\ref{p2.1} for a left-most
anti-oriented infinite
 open path from $(0,0)$, it suffices to prove that, on $\ell^{\rm anti}_{(0,0)}=L^{\rm anti}$,
$|D((0,0),G)|<\infty$ almost surely.

By the proof of Theorem~\ref{th1.1} (i), we know that, on
$\ell^{\rm anti}_{(0,0)}=L^{\rm anti}$, with probability 1, for
any $n\geq 1$, there exists some point $v_n(L^{\rm anti})$ in
$L^{\rm anti}$  from which there is an oriented open path
to $[n,\infty)\times \{0\}$. On the other hand, by (\ref{3.1}) and
the standard ergodic theorem, with probability 1, there exists
infinitely many $m>0$ such that $(m,0)\in\tilde{\fK}$. These
observations imply that, on $\ell^{\rm anti}_{(0,0)}=L^{\rm
anti}$, with probability 1, there exists some $m_0>0$ and $v'\in
L^{\rm anti}$ such that $(m_0,0)\in\tilde{\fK}$ and there is an
oriented open path from $v'$ to some $v''\in
[m_0,\infty)\times\{0\}$. We denote this oriented open path by
$\pi=\pi(v',v'')$.

%Similar to $C_r(\Ga(u,v))$ defined in Section 2, let
%$$C_l(L^{\rm anti}(v(L^{\rm anti}),(0,0))):=\left\{u\in\wedge_{(0,0)}\setminus L^{\rm
%anti}:\ba{ll}& \hskip -5mm u \hb{ lies in the left side of }L^{\rm
%anti}\hb{
% and there is}\\[3mm]
%&\hskip-5mm \hb{a } v\in L^{\rm anti}(v(L^{\rm anti}),(0,0)), \hb{
%such that }u\rightarrow v. \ea\right\},$$ where $L^{\rm
%anti}(v(L^{\rm anti}),(0,0))$ is the finite piece of $L^{\rm
%anti}$ from $v(L^{\rm anti})$ to $(0,0)$.

Let
$$\wedge_{(0,0)}(L^{\rm anti},\pi):=\{u\in\wedge_{(0,0)}\setminus L^{\rm anti}:
u \hb{ lies to the right of } L^{\rm anti} \hb{ and above
 }\pi\}$$ and
$$\bar C_l(L^{\rm anti}(v',(0,0))):=\left\{w\notin
L^{\rm anti}:\ba{ll}&\hskip-5mm w \hb{ lies to the  left of }
L^{\rm anti}\hb{ and for some } z\in L^{\rm
anti}(v',(0,0)),\\[3mm]
&\hskip-5mm z\not=v',\ \ w\rightarrow z\mbox{ uses no open edges
of }L^{\rm anti}\ea\right\}.$$ We declare that
\be\label{3.8}D((0,0),G)\subset \bar C_l(L^{\rm anti}(v',(0,0)))
\cup L^{\rm anti}(v',(0,0))\cup \wedge_{(0,0)}(L^{\rm
anti},\pi).\ee In fact, for any vertex $u$ in $C^{\rm
anti}_{(0,0)}(\subset \fK)$ but outside the set of the right-hand
side of (\ref{3.8}), it is easy to find an open oriented path from
$u$ to $(m_0,0)$. This finding implies $u\notin D((0,0),G)$.

Now, by the definition of left-most anti-oriented infinite open
path, we have \be\label{3.9}|\bar C_l(L^{\rm
anti}(v',(0,0)))|<\infty.\ee Using (\ref{3.8}), (\ref{3.9}), and
the following  fact
$$|L^{\rm anti}(v',(0,0))\cup \wedge_{(0,0)}(L^{\rm anti},\pi)|<\infty,$$ we have
$|D((0,0),G)|<\infty$, so Theorem~\ref{th1.1} (ii) follows. \QED

{\it Proof of Theorem~\ref{th1.1} (iii).} For any $u=(x,y)\in\fK$,
by the standard ergodic theorem,  there exists some
$u'=(x',y')\in \fK$ such that $x'>x$ and $y'=y$ almost surely. It follows  from the
proof of Theorem~\ref{th1.1} (i) that  $\gamma_u$
will meet $\gamma_{u'}$ at some point $v$ of $\fK$ almost surely. Thus, $v$ has
two daughters such that the older one is just the ancestor of $u$ and the other one is her
younger sister.\QED

\section{Proof of Theorem~\ref{th1.4}.}
\renewcommand{\theequation}{4.\arabic{equation}}
\setcounter{equation}{0}

\begin{figure}
\unitlength=0.5mm
\begin{picture}(300,100)(0,-15)
\put(80,0){\vector(1,2){20}} \put(100,40){\vector(2,3){30}}
\put(130,85){\vector(1,3){10}}
\put(80,0){\vector(-2,3){20}}\put(60,30){\vector(-1,3){20}}\put(40,90){\vector(-1,4){8}}
\put(85,-2){$u_1$}

\multiput(79,1)(2,0){1}{\line(1,0){0.5}}\multiput(78,3)(2,0){2}{\line(1,0){0.5}}
\multiput(77,5)(2,0){3}{\line(1,0){0.5}}\multiput(76,7)(2,0){4}{\line(1,0){0.5}}
\multiput(75,9)(2,0){5}{\line(1,0){0.5}}\multiput(74,11)(2,0){6}{\line(1,0){0.5}}
\multiput(72,13)(2,0){7}{\line(1,0){0.5}}\multiput(71,15)(2,0){8}{\line(1,0){0.5}}
\multiput(70,17)(2,0){9}{\line(1,0){0.5}}\multiput(69,19)(2,0){10}{\line(1,0){0.5}}
\multiput(67,21)(2,0){12}{\line(1,0){0.5}}\multiput(66,23)(2,0){13}{\line(1,0){0.5}}
\multiput(64,25)(2,0){14}{\line(1,0){0.5}}\multiput(63,27)(2,0){15}{\line(1,0){0.5}}
\multiput(61,29)(2,0){5}{\line(1,0){0.5}}\multiput(89,29)(2,0){3}{\line(1,0){0.5}}
\multiput(61,31)(2,0){5}{\line(1,0){0.5}}\multiput(89,31)(2,0){4}{\line(1,0){0.5}}
\multiput(59,33)(2,0){6}{\line(1,0){0.5}}\multiput(89,33)(2,0){4}{\line(1,0){0.5}}
\multiput(59,35)(2,0){6}{\line(1,0){0.5}}\multiput(89,35)(2,0){5}{\line(1,0){0.5}}
\multiput(59,37)(2,0){20}{\line(1,0){0.5}}\multiput(59,39)(2,0){20}{\line(1,0){0.5}}
\multiput(57,41)(2,0){22}{\line(1,0){0.5}}
\multiput(57,43)(2,0){15}{\line(1,0){0.5}}\multiput(89,43)(2,0){7}{\line(1,0){0.5}}
\multiput(55,45)(2,0){15}{\line(1,0){0.5}}\multiput(91,45)(2,0){7}{\line(1,0){0.5}}
\multiput(55,47)(2,0){14}{\line(1,0){0.5}}\multiput(93,47)(2,0){6}{\line(1,0){0.5}}
\multiput(55,49)(2,0){13}{\line(1,0){0.5}}\multiput(95,49)(2,0){6}{\line(1,0){0.5}}
\multiput(53,51)(2,0){13}{\line(1,0){0.5}}\multiput(97,51)(2,0){5}{\line(1,0){0.5}}
\multiput(53,53)(2,0){12}{\line(1,0){0.5}}\multiput(99,53)(2,0){5}{\line(1,0){0.5}}
\multiput(53,55)(2,0){11}{\line(1,0){0.5}}\multiput(101,55)(2,0){4}{\line(1,0){0.5}}
\multiput(51,57)(2,0){11}{\line(1,0){0.5}}\multiput(103,57)(2,0){4}{\line(1,0){0.5}}
\multiput(51,59)(2,0){10}{\line(1,0){0.5}}\multiput(105,59)(2,0){4}{\line(1,0){0.5}}
\multiput(51,61)(2,0){9}{\line(1,0){0.5}}\multiput(107,61)(2,0){3}{\line(1,0){0.5}}
\multiput(49,63)(2,0){9}{\line(1,0){0.5}}\multiput(109,63)(2,0){3}{\line(1,0){0.5}}
\multiput(49,65)(2,0){8}{\line(1,0){0.5}}\multiput(111,65)(2,0){3}{\line(1,0){0.5}}
\multiput(49,67)(2,0){7}{\line(1,0){0.5}}\multiput(113,67)(2,0){3}{\line(1,0){0.5}}
\multiput(47,69)(2,0){7}{\line(1,0){0.5}}\multiput(115,69)(2,0){2}{\line(1,0){0.5}}
\multiput(47,71)(2,0){6}{\line(1,0){0.5}}\multiput(117,71)(2,0){2}{\line(1,0){0.5}}
\multiput(47,73)(2,0){5}{\line(1,0){0.5}}\multiput(119,73)(2,0){2}{\line(1,0){0.5}}
\multiput(45,75)(2,0){5}{\line(1,0){0.5}}\multiput(121,75)(2,0){2}{\line(1,0){0.5}}
\multiput(45,77)(2,0){4}{\line(1,0){0.5}}\multiput(123,77)(2,0){1}{\line(1,0){0.5}}
\multiput(45,79)(2,0){3}{\line(1,0){0.5}}\multiput(125,79)(2,0){1}{\line(1,0){0.5}}
\multiput(43,81)(2,0){3}{\line(1,0){0.5}}\multiput(43,83)(2,0){2}{\line(1,0){0.5}}
\multiput(41,85)(2,0){2}{\line(1,0){0.5}}

\put(80,0){\circle*{3}}\put(70,30){$B_{u_1,u_2}$}
\put(85,10){\vector(1,1){20}}\put(97,22){\vector(1,2){5}}%\put(101,30){\vector(1,1){10}}
\put(130,85){\line(-1,-1){44}}\put(87,42){\circle*{3}}\put(83,48){$u_2$}\put(40,89){\line(1,-1){46}}
\put(133,83){$v^r_{1,2}$}\put(130,85){\circle*{2}}
\put(28,85){$v^l_{1,2}$}\put(40,90){\circle*{2}}
\put(55,18){$\ell_{u_1}$}
\put(57,40){\vector(-1,1){20}}\put(45,52){\vector(-1,3){5}}\put(42.5,59.5){\vector(-1,1){8}}
\put(75,-18){$(a)$}\put(115,55){$\gamma_{u_1}$}\put(218,-18){$(b)$}
\put(66,64){$\ell_{u_2}$} \put(100,66){$\gamma_{u_2}$}
\put(200,10){\circle*{3}}\put(198,3){$u_1$}
\multiput(199,13)(2,0){2}{\line(1,0){0.5}}\multiput(199,15)(2,0){2}{\line(1,0){0.5}}
\multiput(197,17)(2,0){4}{\line(1,0){0.5}}\multiput(197,19)(2,0){4}{\line(1,0){0.5}}
\multiput(195,21)(2,0){5}{\line(1,0){0.5}}\multiput(195,23)(2,0){6}{\line(1,0){0.5}}
\multiput(193,25)(2,0){7}{\line(1,0){0.5}}\multiput(193,27)(2,0){8}{\line(1,0){0.5}}
\multiput(191,29)(2,0){9}{\line(1,0){0.5}}\multiput(191,31)(2,0){10}{\line(1,0){0.5}}
\multiput(189,33)(2,0){11}{\line(1,0){0.5}}\multiput(189,35)(2,0){12}{\line(1,0){0.5}}
\multiput(189,37)(2,0){13}{\line(1,0){0.5}}\multiput(187,39)(2,0){14}{\line(1,0){0.5}}
\multiput(187,41)(2,0){14}{\line(1,0){0.5}}\multiput(187,43)(2,0){14}{\line(1,0){0.5}}
\multiput(185,45)(2,0){14}{\line(1,0){0.5}}\multiput(185,47)(2,0){14}{\line(1,0){0.5}}
\multiput(185,49)(2,0){1}{\line(1,0){0.5}}\multiput(205,49)(2,0){3}{\line(1,0){0.5}}
\multiput(183,51)(2,0){2}{\line(1,0){0.5}}\multiput(205,51)(2,0){3}{\line(1,0){0.5}}
\multiput(183,53)(2,0){2}{\line(1,0){0.5}}\multiput(205,53)(2,0){2}{\line(1,0){0.5}}
\multiput(183,55)(2,0){2}{\line(1,0){0.5}}\multiput(205,55)(2,0){1}{\line(1,0){0.5}}
\multiput(181,57)(2,0){13}{\line(1,0){0.5}}\multiput(181,59)(2,0){12}{\line(1,0){0.5}}
\multiput(181,61)(2,0){12}{\line(1,0){0.5}}\multiput(181,63)(2,0){11}{\line(1,0){0.5}}
\multiput(181,65)(2,0){11}{\line(1,0){0.5}}\multiput(181,67)(2,0){10}{\line(1,0){0.5}}
\multiput(181,69)(2,0){10}{\line(1,0){0.5}}\multiput(181,71)(2,0){10}{\line(1,0){0.5}}
\multiput(181,73)(2,0){9}{\line(1,0){0.5}}\multiput(181,75)(2,0){9}{\line(1,0){0.5}}
\multiput(181,77)(2,0){8}{\line(1,0){0.5}}\multiput(181,79)(2,0){8}{\line(1,0){0.5}}
\multiput(181,81)(2,0){7}{\line(1,0){0.5}}\multiput(181,83)(2,0){7}{\line(1,0){0.5}}
\multiput(181,85)(2,0){6}{\line(1,0){0.5}}\multiput(181,87)(2,0){6}{\line(1,0){0.5}}
\multiput(181,89)(2,0){5}{\line(1,0){0.5}}\multiput(181,91)(2,0){5}{\line(1,0){0.5}}
\multiput(181,93)(2,0){4}{\line(1,0){0.5}}\multiput(181,95)(2,0){4}{\line(1,0){0.5}}
\multiput(181,97)(2,0){3}{\line(1,0){0.5}}\multiput(181,99)(2,0){3}{\line(1,0){0.5}}
\multiput(181,101)(2,0){2}{\line(1,0){0.5}}\multiput(181,103)(2,0){2}{\line(1,0){0.5}}
\multiput(181,103)(2,0){1}{\line(1,0){0.5}}\multiput(181,105)(2,0){1}{\line(1,0){0.5}}

\put(202.5,15){\vector(1,1){12}}\put(208,20){\vector(1,3){3}}\put(209,24){\vector(1,1){8}}
\put(250,0){\circle*{3}}\put(254,0){$u_2$}
\put(200,10){\vector(-1,2){10}}\put(190,30){\vector(-1,3){10}}\put(180,60){\vector(0,1){40}}
\put(180,100){\vector(0,1){25}}
\put(200,10){\vector(1,2){15}}\put(215,40){\vector(2,3){20}}
\put(235,70){\vector(1,1){40}}\put(250,0){\vector(0,1){40}}\put(250,40){\vector(1,4){15}}
\put(250,0){\vector(-1,1){20}}\put(230,20){\vector(-2,3){20}}\put(210,50){\vector(-1,2){30}}
\put(250,100){$v^r_{1,2}$}\put(265,100){\circle*{2}}
\put(165,110){$v^l_{1,2}$}\put(180,110){\circle*{2}}
\put(216,41){\circle*{2}}\put(220,40){$v_{1,2}$}
\put(220,65){$\gamma_{u_1}$} \put(254,50){$\gamma_{u_2}$}
\put(178,33){$\ell_{u_1}$} \put(180,60){\vector(-1,1){15}}
\put(140,50){\small $\mbox{a left bud}$}
\put(150,57){\vector(1,1){8}}\put(168,74){\circle{22}}\put(172,68){\vector(0,1){10}}\put(172,73){\vector(-2,3){5}}
\put(186,50){$\Delta_{u_1,u_2}$}\put(200,73){$\ell_{u_2}$}
\put(118,05){\vector(-1,1){8}}\put(105,-1){\small $\mbox{a right
bud}$}\put(104,23){\circle{21}}
\end{picture}
\begin{center}
\begin{minipage}{13.5cm}
\caption{\small $B_{u_1,u_2}$, $\Delta_{u_1,u_2}$, and the left
(resp., right) {\it buds} planted in $\ell_{u_1}$ (resp.,
$\gamma_{u_1}$) are all finite. Note that
$C_l(\ell_{u_1}(u_1,v^l_{1,2}))$ consists of all such finite left
buds planted in $\ell_{u_1}(u_1,v^l_{1,2})\setminus\{v^l_{1,2}\}$,
so it is finite. The situations for
$C_r(\gamma_{u_1}(u_1,v^r_{1,2}))$ and
$C_r(\gamma_{u_1}(u_1,v_{1,2}))$ are the same.} \label{f2}
\end{minipage}
\end{center}
\end{figure}
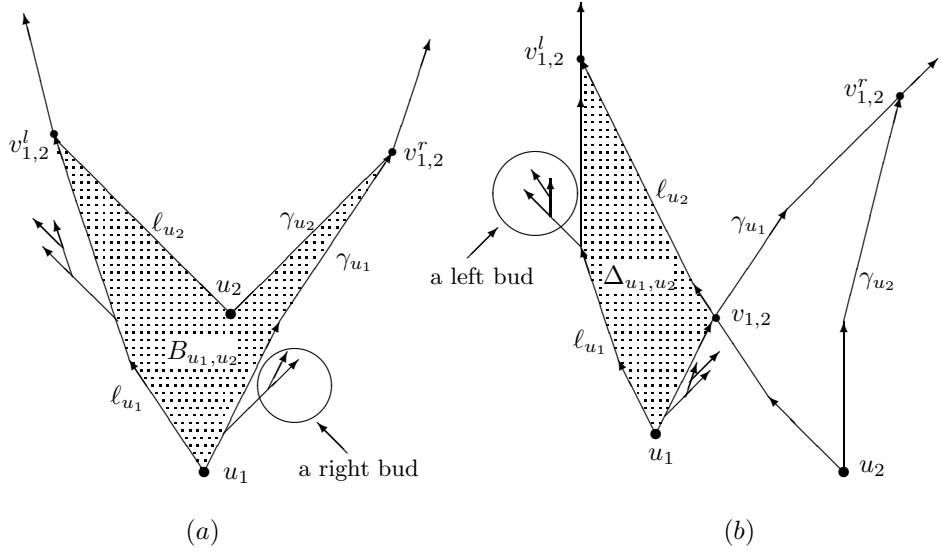
It suffices to prove that, for
any $u_1,u_2\in \fK$, $|C_{u_1}\Delta C_{u_2}|<\infty$, where
$C_{u_1}\Delta C_{u_2}$ is the symmetric difference of $C_{u_1}$
and $C_{u_2}$. By Theorem~\ref{th1.1}, with probability 1,
$$\gamma_{u_1}\cap\gamma_{u_2}\not=\emptyset\mbox{ and } \ell_{u_1}\cap\ell_{u_2}\not=\emptyset.$$
Let
$$v^r_{1,2}=(x_r, y_r)\in
\gamma_{u_1}\cap\gamma_{u_2}\mbox{ and }v^l_{1,2}=(x_l,y_l)\in
\ell_{u_1}\cap\ell_{u_2}$$
 be the vertices with the smallest second coordinates. With these definitions,
we will prove that \be\label{4.1}|C_{u_1}\setminus
C_{u_2}|<\infty.\ee By (\ref{4.1}) and  symmetry, we also have
$|C_{u_1}\setminus C_{u_2}|<\infty$, so Theorem~\ref{th1.4}
follows.

Now it remains to show (\ref{4.1}). Without loss of generality, we
divide the problem into following two cases (see Fig.~\ref{f2}):
\begin{enumerate}
    \item $u_2$ lies within the cone between
$\ell_{u_1}$ and $\gamma_{u_1}$;
    \item $u_2$ does not lie within the
cone between $\ell_{u_1}$ and $\gamma_{u_1}$.
\end{enumerate}

We focus on case 1 (see Fig.~\ref{f2} (a)). Let $B_{u_1,u_2}$ be
the finite {\it butterfly} shape enclosed by
$\gamma_{u_2}(u_2,v^r_{1,2})$, $\gamma_{u_1}(u_1,v^r_{1,2})$,
$\ell_{u_2}(u_2,v^l_{1,2})$, $\ell_{u_1}(u_1,v^l_{1,2})$, and the
vertices surrounded by them.  It is clear that
$$ C_{u_1}\setminus
C_{u_2}\subset B_{u_1,u_2}\cup
C_r(\gamma_{u_1}(u_1,v^r_{1,2}))\cup
C_l(\ell_{u_1}(u_1,v^l_{1,2})).$$ Note that
$C_r(\gamma_{u_1}(u_1,v^r_{1,2}))$ and
$C_l(\ell_{u_1}(u_1,v^l_{1,2}))$ are defined in (\ref{cr}). By the
definition of $\gamma_{u_1}$ and $\ell_{u_1}$, we have
$$|C_r(\gamma_{u_1}(u_1,v^r_{1,2}))|<\infty;\ \ |C_l(\ell_{u_1}(u_1,v^l_{1,2}))|<\infty.$$
This tells us that $|C_{u_1}\setminus C_{u_2}|<\infty$, so
(\ref{4.1}) follows when case 1 holds.

Let us focus on case 2 (see Fig.~\ref{f2} (b)). Without loss of
generality, we may further assume that $u_1$ and $u_2$ have the
relative position such that $\gamma_{u_1}\cap
\ell_{u_2}\not=\emptyset$. Let $v_{1,2}\in \gamma_{u_1}\cap
\ell_{u_2}$ be the vertex with the smallest second coordinate.
Moreover, let $\Delta_{u_1,u_2}$ be the finite {\it triangle}
shape enclosed by  $\gamma_{u_1}(u_1,v_{1,2})$,
$\ell_{u_1}(u_1,v^l_{1,2})$, $\ell_{u_1}(v_{1,2},v^l_{1,2})$, and
the vertices surrounded by them. It is clear that
$$ C_{u_1}\setminus
C_{u_2}\subset \Delta(u_1,u_2)\cup
C_r(\gamma_{u_1}(u_1,v_{1,2}))\cup
C_l(\ell_{u_1}(u_1,v^l_{1,2})).$$ The same argument for the first
case tells us  that $|C_{u_1}\setminus C_{u_2}|<\infty$, so
(\ref{4.1}) also follows when case 2 holds.\QED

\vskip 5mm
\begin{minipage}{7cm}
{\noindent\hskip-3mm $^{~1}$Department of Mathematics, Capital
Normal University, 100037, Beijing, China. E-mail:
\texttt{wuxy@mail.cnu.edu.cn} \vskip 2mm
\noindent\hskip-3mm$^{~2}$ Department of Mathematics, University
of Colorado Colorado Springs, CO 80933-7150. E-mail:
\texttt{yzhang@math.uccs.edu}}
\end{minipage}
\end{document}